\documentstyle[11pt]{article}
\title{$\eta$-invariant and Chern-Simons current}
\author{Weiping Zhang\thanks{Partially
supported by  the MOEC  and the 973 project.}}
\date{  }
\begin{document}
\maketitle

\begin{abstract}
We show that the  ${\bf R}/{\bf Z}$ part of the analytically
defined  $\eta$ invariant of Atiyah-Patodi-Singer for a Dirac
operator on an odd dimensional closed spin manifold  can be
expressed purely geometrically through a
 stable Chern-Simons current on a higher dimensional sphere. As a
 preliminary application, we discuss the relation with the
 Atiyah-Patodi-Singer ${\bf R}/{\bf Z}$
 index theorem for unitary flat vector  bundles, and prove an ${\bf R}$ refinement
 in the case where the Dirac operator is replaced by the Signature operator. We also extend the
 above discussion to the case of $\eta$ invariants associated to
 Hermitian vector bundles with non-unitary connection, which are
 constructed by using a trick due to Lott.
\end{abstract}

$\ $

\noindent{\bf\S 1. Introduction}

$\ $

The $\eta$ invariant of Atiyah-Patodi-Singer was introduced in
[APS1] as the correction term on the boundary of the index theorem
for Dirac operators on manifolds with boundary. Since then it has
appeared in many parts of geometry, topology as well as physics.
We first recall the definition of this $\eta$ invariant.

Let $M$ be an odd dimensional closed oriented spin Riemannian
manifold. Let $S(TM)$ be  the associated bundle of spinors. Let
$E$ be a Hermitian vector bundle over $M$ carrying with a
Hermitian connection. Then one can define canonically a Dirac
operator $D^E:\Gamma(S(TM)\otimes E) \rightarrow
\Gamma(S(TM)\otimes E)$. It is a formally self-adjoint first
order elliptic differential operator.

Let $s\in {\bf C}$ with ${\rm Re}(s)>{\dim M\over 2}$. Following
[APS1], one defines the $\eta$ function of $D^E$ by
$$\eta(D^E,s)=\sum_{\lambda\in {\rm Spec}(D^E)\setminus \{ 0 \}} {{\rm sgn}(\lambda)\over
|\lambda|^s}. \eqno(1.1)$$ It is shown in [APS1] that
$\eta(D^E,s)$ is a holomorphic function for ${\rm Re}(s)>{\dim
M\over 2}$, and can be extended to a meromorphic function on
${\bf C}$. Moreover, it is holomorphic at $s=0$. The value of
$\eta(D^E,s)$ at $s=0$ is called the $\eta$ invariant of $D^E$
and is denoted by $\eta(D^E)$. Let $\bar{\eta}(D^E)$ be  the
reduced $\eta$ invariant of $D^E$ which is also defined in [APS1]:
$$\bar{\eta}(D^E)={\dim (\ker D^E)+\eta(D^E)\over
2}.\eqno(1.2)$$

It turns out that this analytically defined invariant may jump by
integers as the metrics and connections on $TM$ and $E$ change.
These jumps can be detected by  spectral flows introduced in
[APS3]. On the other hand, the mod ${\bf Z}$ component of
$\bar{\eta}(D^E)$ is smooth with respect to the involved metrics
and connections, and its variation can be expressed through
Chern-Simons forms. However, whether $\bar{\eta}(D^E)$ mod ${\bf
Z}$ itself can be expressed geometrically, without passing to the
spectral set of $D^E$, remains a question for long time. Here we
only mention that such a formula for $\bar{\eta}(D^E)$ mod ${\bf
Q}$ was proved
 in [CS, Theorem 9.1], expressing the ${\bf R}/{\bf Q}$ component
of $\bar{\eta}(D^E)$ through the Cheeger-Simons differential
characters.

The purpose of this short article  is to show that there is indeed
a purely geometric formula for $\bar{\eta}(D^E)$ mod ${\bf Z}$.
More precisely, if we embed  $M$ into a higher odd dimensional
sphere, then $\bar{\eta}(D^E)$ mod ${\bf Z}$ can be expressed
through a Chern-Simons current on the sphere. Comparing with
Cheeger-Simons' mod ${\bf Q}$ one, such a formula is  more of
$K$-theoretic nature, and should be viewed as an index theorem in
some geometric $K$-theory (compare with [AS] and [L]).

In fact, this formula, which will be stated in its precise form
in Theorem 2.3, can be obtained as an immediate application of a
localization formula  for $\eta$ invariants proved by Bismut and
Zhang in [BZ]. Our simple observation is that if one  applies the
Bismut-Zhang formula to an embedding into a higher dimensional
sphere, then a simple application of the Bott periodicity will
lead us to a geometric formula for the ${\bf R}/{\bf Z}$
component of $\bar{\eta}(D^E)$.

Thus in the next section we will recall the Bismut-Zhang
localization formula for $\eta$ invariants and prove the
geometric formula for $\bar{\eta}(D^E)$ mod ${\bf Z}$.

Also recall that the proof given by Bismut-Zhang in [BZ] for their
localization formula relies heavily on the difficult paper of
Bismut-Lebeau [BL], and might not be easy to follow. So in
Section 3 we will give an alternate  proof of this localization
formula by making use of the Freed-Melrose index theorem for
${\bf Z}/k$ manifolds [FM] instead.

In Section 4, we present some preliminary applications of our
formula to the case of flat vector bundles. In particular, we
show that our formula leads to an intrinsic formulation of the
Atiyah-Patodi-Singer ${\bf R}/{\bf Z}$ index theorem for unitary
flat vector bundles [APS3, Theorem 5.3]. Moreover, we show that
when considering the Signature operator instead of the Dirac
operator,  one can refine the above index theorem to an ${\bf R
}$ valued one.


Now note that it is pointed out in [APS3] that the above mentioned
index theorem for flat vector bundles [APS3, Theorem 5.3] indeed
holds for all flat connections not necessarily unitary, and that
one can extend the definition of the $\eta$ invariant for these
not necessarily unitary connections.

In Section 5 of this paper, we will propose a definition of the
$\eta$ invaraint for Dirac operators coupled with Hermitian
vector bundles with non-unitary connection, by using a trick of
Lott in [L, Section 5]. When applying to flat connections, we
show that this invariant verifies the Atiyah-Patodi-Singer index
theorem [APS3, Theorem 5.3], which implies that when mod ${\bf
Z}$, our construction agrees with what indicated in [APS3].  We
will discuss further applications of this invaraint elsewhere.

$$\ $$

\noindent {\bf \S 2. A geometric formula for $\eta$ invariants}

$\ $

In this section, we recall the localization formula for $\eta$
invariants of Bismut-Zhang [BZ] and use it to deduce a geometric
formula for  $\eta$ invariants.

This section is organized as follows. In a), we recall the direct
image construction of Atiyah-Hirzebruch [AH] under real embeddings
in a geometrical form. In b),  we recall from [B2] and [BZ] the
construction of the Chern-Simons current associated to the
geometric direct image constructed in a). In c), we state the
localization formula from [BZ]. In d), we apply this localization
formula to get a geometric formula expressing the $\eta$
invariants through Chern-Simons currents on spheres.

$\ $

{\bf a). A geometric construction of direct images}

$\ $

Let $i:Y\hookrightarrow X$ be an embedding between two smooth
oriented manifolds. We make the assumption that $\dim X-\dim Y$
is even and that if $N$ denotes the normal bundle to $Y$ in $X$,
then $N$    is orientable,  spin and carries an induced
orientation as well as a (fixed) spin structure.

 Let  $\mu$ be a complex vector bundle over $Y$.

  Atiyah
and Hirzebruch have constructed  in [AH] an element $i_!\mu\in
\widetilde{K}(X)$,  called the direct image of $\mu$ under $i$. We
here recall this construction in a geometric form.

Let $g^N$ be a Euclidean metric on $N$ and $\nabla^N$  a
Euclidean connection on $N$ preserving $g^N$. Let $S(N)$ be the
vector bundle of spinors associated to $(N,g^N)$. Then
$S(N)=S_+(N)\oplus S_-(N)$ (resp. its dual $S^*(N)=S^*_+(N)\oplus
S^*_-(N)$) is a ${\bf Z}_2$-graded complex vector bundle over $Y$
carrying with an induced Hermitian metric
$g^{S(N)}=g^{S_+(N)}\oplus g^{S_-(N)}$ (resp.
$g^{S^*(N)}=g^{S^*_+(N)}\oplus g^{S^*_-(N)}$) from $g^N$, as well
as a Hermitian connection
$\nabla^{S(N)}=\nabla^{S_+(N)}\oplus\nabla^{S_-(N)}$ (resp.
$\nabla^{S^*(N)}=\nabla^{S^*_+(N)}\oplus\nabla^{S^*_-(N)}$)
induced from $\nabla^N$.

Let $g^\mu$ be a Hermitian metric on $\mu$ and $\nabla^\mu$  a
Hermitian connection on $\mu$ preserving $g^\mu$.

For any $r>0$, set $N_r=\{Z\in N: |Z|< r\}.$ We make the
assumption that there is $\varepsilon_0>0$ such that
$N_{2\varepsilon_0}$ is diffeomorphic to an open neighborhood
 of $Y$ in $X$. Without confusion we now view directly
 $N_{2\varepsilon_0}$ as an open neighborhood of $Y$ in $X$.

 Let $\pi:N\rightarrow Y$ denote the projection of the normal
 bundle $N$ over $Y$.

 If $Z\in  {N}$, let $\tilde{c}(Z)\in {\rm End}(S^*( {N})) $
 be the transpose of $c(Z)$
acting on $S( {N})$. Let $\tau^{ {N}*}\in {\rm End}(S^*( {N})) $
be the transpose of $\tau^{ {N}}$ defining the
 ${\bf Z}_2$-grading of $S( {N})=S_+( {N})\oplus S_-( {N})$.

 Let $\pi^*(S^*(N))$ be the pull back bundle of $S^*(N)$ over $N$.
For any $Z\in N$ with $Z\neq 0$, let $\tau^{ {N}*}\tilde{c}(Z):
 \pi^*(S^*_\pm(N))|_Z\rightarrow \pi^*(S^*_\mp(N))|_Z$ denote the
 corresponding pull back isomorphisms at $Z$.

Let $F$ be a complex vector bundle over $Y$ such that
$S_-(N)\otimes\mu\oplus F$ is a  trivial complex vector bundle
over $Y$ (cf. [A]). Then
$$\tau^{ {N}*}\tilde{c}(Z)\oplus\pi^*{\rm
Id}_F :\pi^*(S_+^*(N)\otimes\mu\oplus F)\rightarrow
\pi^*(S_-^*(N)\otimes\mu\oplus F)\eqno(2.1)$$ induces an
isomorphism between two trivial  vector bundles over $
N_{2\varepsilon_0}\setminus Y$.

Let $F$ admit  a Hermitian metric $g^F$ and  a  Hermitian
connection $\nabla^F$.

Clearly, $\pi^*(S_\pm^* (N)\otimes\mu\oplus F)|_{\partial
N_{2\varepsilon_0}}$ extend smoothly to two trivial complex vector
bundles over $X\setminus N_{2\varepsilon_0}$. Moreover, the
isomorphism $\tau^{ {N}*}\tilde{c}(Z)\oplus\pi^*{\rm Id}_F$ over
$\partial N_{2\varepsilon_0}$ extends smoothly to an isomorphism
between these two trivial vector bundles over $X\setminus
N_{2\varepsilon_0}$.

In summary, what we get is a ${\bf Z}_2$-graded Hermitian vector
bundle $$\xi=\xi_+\oplus \xi_-,\ \ \ \ g^\xi=g^{\xi_+}\oplus
g^{\xi_-}\eqno(2.2)$$ over $X$ such that
$$\xi_\pm|_{N_{\varepsilon_0}} =\pi^*(S_\pm^*(N)\otimes\mu\oplus
F)|_{N_{\varepsilon_0}},\ \ \ \ g^{\xi_\pm|_{N_{\varepsilon_0}}}
=\pi^*(g^{S_\pm^*(N)\otimes\mu}\oplus
g^F)|_{N_{\varepsilon_0}},\eqno(2.3)$$ where
$g^{S_\pm^*(N)\otimes\mu}$ is the tensor product Hermitian metric
on $S_\pm^*(N)\otimes\mu$ induced from $g^{S_{\pm}^*(N)}$ and
$g^\mu$. It is easy to see that there exists an odd self-adjoint
automorphism $V$ of $\xi$ such that
$$V|_{N_{\varepsilon_0}} =\tau^{ {N}*}\tilde{c}(Z)\oplus\pi^*{\rm
Id}_F.\eqno(2.4)$$ Moreover, there is a ${\bf Z}_2$-graded
Hermitian connection $\nabla^\xi=\nabla^{\xi_+}\oplus
\nabla^{\xi_-}$ on $\xi=\xi_+\oplus\xi_-$ over $X$ such that
$$\nabla^{\xi_\pm}|_{N_{\varepsilon_0}}=\pi^*(\nabla^{S^*_\pm(N)\otimes\mu}\oplus\nabla^F),\eqno(2.5)$$
where $\nabla^{S^*_\pm(N)\otimes\mu}$ is the Hermitian connection
on $\nabla^{S^*_\pm(N)\otimes\mu}$ defined by
$\nabla^{S^*_\pm(N)\otimes\mu}=\nabla^{S^*_\pm(N)}\otimes {\rm
Id}_\mu +{\rm Id}_{S^*_\pm (N)}\otimes\nabla^\mu$.

Clearly, $\xi_+-\xi_-\in \widetilde{K}(X)$ is exactly the
Atiyah-Hirzebruch direct image $i_!\mu$ of $\mu$ constructed in
[AH]. We call $(\xi,\nabla^\xi, V)$ constructed above a geometric
direct image of ($\mu$, $\nabla^\mu$).

$\ $

{\bf b). A Chern-Simons current associated to a geometric direct
image }

$\ $

We make the same assumptions and use the same notations as in a).

If $E$ is a real vector bundle carrying with a connection
$\nabla^E$, we denote by $\widehat{A}(E,\nabla^E)$ the Hirzebruch
characteristic form defined by
$$\widehat{A}(E, \nabla^{E})
={\det}^{1/2}\left({{\sqrt{-1}\over  4\pi}R^{E} \over \sinh\left({
\sqrt{-1}\over 4\pi}R^{E}\right)}\right),\eqno(2.6)$$ where
$R^E=\nabla^{E,2}$ is the curvature of $\nabla^E$. While if $E'$
is a complex vector bundle carrying with a connection
$\nabla^{E'}$, we denote by ${\rm ch}(E',\nabla^{E'})$ the Chern
character form associated to $(E',\nabla^{E'})$ (cf. [Z2, Section
1]).

Let $i^{1/2}$ be a fixed square root of $i=\sqrt{-1}$. The
objects which will be considered in the sequel do not depend on
this square root. Let $\varphi$ be the map $\alpha\in
\Lambda^*(T^*X)\rightarrow (2\pi i)^{-{\deg \alpha\over 2}}\alpha
\in\Lambda^*(T^*X)$.

We now use Quillen's superconnection formalism. For $T\geq 0$, let
$C_T$ be the superconnection on the ${\bf Z}_2$-graded vector
bundle $\xi$ defined by
$$C_T=\nabla^\xi+\sqrt{T}V.\eqno (2.7)$$
The curvature $C_T^2$ of $C_T$ is a smooth section of
$(\Lambda^*(T^*M)\widehat{\otimes}{\rm End}(\xi))^{\rm even}$. By
[Q], we know that for any $T>0$, $${\partial\over\partial T}{\rm
Tr}_s\left[\exp\left(-C^2_T\right)\right]=-{d\over 2\sqrt{T}}{\rm
Tr}_s\left[V\exp\left(-C^2_T\right)\right].\eqno(2.8)$$

Clearly, the technical assumptions in [BZ, (1.10)-(1.12)] hold for
our constructions in a). Thus one can proceed as in [B1], [B2] and
[BZ, Definition 1.3] to construct the Chern-Simons current
$\gamma^{\xi, V}$ as
$$\gamma^{\xi,
V}={1\over \sqrt{2\pi i}}\int_0^{+\infty}\varphi{\rm
Tr}_s\left[V\exp\left(-C^2_T\right)\right]{dT\over
2\sqrt{T}}.\eqno(2.9)$$

Let $\delta_Y$ denote the current of integration over the
oriented submanifold $Y$ of $X$. Then by [BZ, Theorem 1.4], we
have that
$$d\gamma^{\xi,
V}={\rm ch}(\xi_+,\nabla^{\xi_+})-{\rm
ch}(\xi_-,\nabla^{\xi_-})-{\widehat{A}}^{-1}(N,\nabla^N){\rm ch
}(\mu,\nabla^\mu)\delta_Y.\eqno(2.10)$$ Moreover, as indicated in
[BZ, Remark 1.5], by proceeding as in [BGS, Theorem 3.3], one can
prove that $\gamma^{\xi, V}$ is a locally integrable current.

$\ $

 {\bf c). A localization formula for $\eta$ invariants}

 $\ $

 We assume in this subsection that $i:Y\hookrightarrow X$ is an
 embedding between two odd dimensional closed oriented spin manifolds.
 Then the normal bundle $N$ to $Y$ in $X$ is  even
 dimensional and carries a canonically induced orientation and spin
 structure.
 Let $g^{TX}$  be a Riemannian metric on $TX$.
 Let $g^{TY}$ be the restricted Riemannian metric on
$TY$. Let $\nabla^{TX}$ (resp. $\nabla^{TY}$) denote the
 Levi-Civita connection associated to $g^{TX}$ (resp. ${g^{TY}}$).

Without loss of generality we may and we will make the assumption
that the embedding $(Y,g^{TY})\hookrightarrow (X,g^{TX})$ is
totally geodesic. Let $N$ carry the canonically induced Euclidean
metric as well as the  Euclidean connection.

Thus,  we may and we will make the same construction as in a), b).

Recall that the definition of the reduced $\eta$ invariant for a
(twisted) Dirac operator on an odd dimensional spin Riemannian
manifold has been recalled in Section 1.

Under our assumptions, we see easily that the localization
formula for $\eta$ invariants proved in [BZ] holds in a slightly
simplified form. We recall it as follows.

$\ $

\noindent {\bf Theorem 2.1} (Bismut-Zhang [BZ, Theorem 2.2]) {\it
The following identity holds, $$\bar{\eta}(D^{\xi_+}) -
\bar{\eta}(D^{\xi_-})\equiv\bar{\eta}(D^\mu)+\int_X\widehat{A}(TX,\nabla^{TX})\gamma^{\xi,V}\
\ \ {\rm mod}\ {\bf Z}.\eqno(2.11)$$ }


\noindent {\bf Remark 2.2} The extra Chern-Simons form in [BZ,
Theorem 2.2] disappears here simply because we have made the
simplifying assumption that the isometric embedding
$(Y,g^{TY})\hookrightarrow (X,g^{TX})$ is totally geodesic.

$\ $

{\bf d). A geometric formula for $\bar{\eta}(D^\mu)$}

$\ $

We continue the discussion in c) and assume that $X=S^{2n-1}$, a
higher odd dimensional sphere (but we do not assume that it
admits the standard metric, this makes the isometric embedding
$i:(Y,g^{TY})\hookrightarrow (S^{2n-1},g^{TS^{2n-1}})$ to be
totally geodesic possible).

Now recall that by the Bott periodicity (cf. [1]), one has
$\widetilde{K}(S^{2n-1})=\{ 0 \}$. Thus, in our case, we have
$i_!\mu=0$. This means there is a trivial complex vector bundle
$\theta$ over $S^{2n-1}$ such that $\theta\oplus \xi_+$ is
isomorphic to $\theta\oplus \xi_-$.

We equip $\theta$ with a Hermitian metric as well as a Hermitian
connection $\nabla^\theta$.

Let ${\xi}'= {\xi}_+'\oplus{\xi}_-'$ be the ${\bf Z}_2$-graded
Hermitian vector bundle over $S^{2n-1}$ with
${\xi}_\pm'=\theta\oplus \xi_\pm$. Then $\xi'_\pm$ and $\xi'$
carry canonically induced Hermitian connections
$\nabla^{\xi'_\pm}$ and  $\nabla^{\xi'}$ respectively  through
direct sums.

Let $W$ be an odd self-adjoint automorphism of ${\xi}'$, which
clearly exists by the above discussion.

For any $T\geq 0$, set
$$C'_T=\nabla^{\xi'}+\sqrt{T}W.\eqno(2.12)$$

Similarly as in (2.9), we  define the Chern-Simons form
$\gamma^{\xi',W}$ as
$$\gamma^{\xi',W} ={1\over \sqrt{2\pi i}}\int_0^{+\infty}\varphi{\rm
Tr}_s\left[W\exp\left(-\left(C'_T\right)^2\right)\right]{dT\over
2\sqrt{T}}.\eqno(2.13)$$

Since $W$ is invertible, one has the following formula due to
Bismut-Cheeger  [BC, Theorem 2.28], which corresponds to the case
with $Y=\emptyset$ in (2.11),
$$\bar{\eta}(D^{\xi'_+})-\bar{\eta}(D^{\xi'_-})\equiv \int_{S^{2n-1}}
\widehat{A}(TS^{2n-1},\nabla^{TS^{2n-1}})\gamma^{\xi',W}\ \ \ {\rm
mod}\ {\bf Z}.\eqno(2.14)$$ On the other hand, one clearly has
that
$$\bar{\eta}(D^{\xi'_\pm})=\bar{\eta}(D^{\xi_\pm})+\bar{\eta}(D^\theta).\eqno(2.15)$$

{}From (2.11), (2.14) and (2.15), one gets
$$\bar{\eta}(D^\mu)\equiv\int_{S^{2n-1}}
\widehat{A}(TS^{2n-1},\nabla^{TS^{2n-1}})\gamma^{\xi',W}
-\int_{S^{2n-1}}\widehat{A}(TS^{2n-1},\nabla^{TS^{2n-1}})\gamma^{\xi,V}\
\ \ {\rm mod}\ {\bf Z}.\eqno(2.16)$$

We can re-formulate (2.16) as follows.

Let $\widetilde{\xi}=\widetilde{\xi}_+\oplus\widetilde{\xi}_-$ be
the ${\bf Z}_2$-graded Hermitian vector bundle over $S^{2n-1}$
defined by $$\widetilde{\xi}_+=\xi_+\oplus\xi_-',\ \ \ \
\widetilde{\xi}_-=\xi_-\oplus\xi_+',\eqno(2.17)$$ carrying the
canonically induced Hermitian connection
$\nabla^{\widetilde{\xi}}$ through direct sums. Let
$\widetilde{V}=V\oplus W^T$, where $W^T$ is the transpose of $W$,
be the odd self-adjoint automorphism of $\widetilde{\xi}$. Then
$(\widetilde{\xi},\nabla^{\widetilde{\xi}},\widetilde{V} )$ forms
a geometric direct image of $\mu$ in the sense of Section 2a).

Let $\gamma^{\widetilde{\xi},\widetilde{V}}$ be the associated
Chern-Simons current defined by (2.9). By (2.10) and the
construction of
$(\widetilde{\xi},\nabla^{\widetilde{\xi}},\widetilde{V} )$, one
verifies easily that
$$d\gamma^{\widetilde{\xi},
\widetilde{V}}=-{\widehat{A}}^{-1}(N,\nabla^N){\rm ch
}(\mu,\nabla^\mu)\delta_Y.\eqno(2.18)$$

We can now state  the main result of this section as follows,
which is simply a re-formulation of (2.16).

$\ $

\noindent {\bf Theorem 2.3} {\it The following identity holds,
$$\bar{\eta}(D^\mu)\equiv -\int_{S^{2n-1}}\widehat{A}(TS^{2n-1},\nabla^{TS^{2n-1}})
\gamma^{\widetilde{\xi},\widetilde{V}}\ \ \ {\rm mod}\ {\bf
Z}.\eqno(2.19)$$}

\noindent {\bf Remark 2.4} It is remarkable that the right hand
side of (2.19) does not involve any spectral information of
$D^\mu$. It is purely geometric/topological, although the
existence of the invertible element $W$ is by no means trivial
(as we have seen, it follows from the Bott periodicity). It
resembles well the $K$-theoretic proof of the Atiyah-Singer index
theorem [AS]. On the other hand, the different choice of $W$ may
cause the right hand side of (2.19) an integer jump. This partly
explains that (2.19) is in general a mod ${\bf Z}$ formula. While
conversely, one can always find a $W$ to make (2.19) a purely
equality in ${\bf R}$. This may sound few sense as when $g^{TM}$,
$g^\mu$ and $\nabla^\mu$ vary, $\dim(\ker D^\mu)$ may jump.
However, if we take $\mu=S(TY)$, then $D^{S(TY)}$ is the Signature
operator of $(Y,g^{TY})$ and such an equality in ${\bf R}$ with a
suitable choice of $W$ will not depend on the variation of
$g^{TY}$.

$\ $

\noindent {\bf Remark 2.5} If $\dim Y\equiv 3\ {\rm mod}\ 8{\bf
Z}$ and $\mu$ is a complexification of a Euclidean vector bundle
carrying with a Euclidean connection, then one can embed $Y$ into
a higher $8l+3$ dimensional manifold and proceed as in [Z1,
Section 3] to improve (2.19) to  a mod $2{\bf Z}$ formula.

$\ $

\noindent {\bf Remark 2.6} As have been mentioned in Section 1,
the proof of Theorem 2.1 given in [BZ] relies heavily on the
difficult paper of Bismut-Lebeau [BL]. So in the next section, we
will give an alternate proof of Theorem 2.1 for the case where
$X$ is a higher dimensional sphere.

$$\ $$

\noindent {\bf \S 3. An alternate proof of Theorem 2.1}

$\ $


As in Section 2, let $Y$ be an odd dimensional closed oriented
spin manifold carrying with a Riemannian metric $g^{TY}$ and the
associated Levi-Civita connection $\nabla^{TY}$. Let $\mu$ be a
complex vector bundle over $Y$ carrying with a Hermitian metric
$g^\mu$ and a Hermitian connection $\nabla^\mu$.

In case when there will be no confusion, we will use the
notations in Section 2 without further explanation.

Since $\dim Y$ is odd, by a well-known result in
bordism/cobordism theory, there is a positive integer $k$ such
that the $k$ disjoint copies of $Y$ bound a compact oriented spin
manifold $\widehat{Y}$ of dimension $\dim Y+1$ such that the
boundary $\partial \widehat{Y}$ does not contain other components.
Moreover, there is a  complex vector bundle $\widehat{\mu}$ over
$\widehat{Y}$ such that when restricted to boundary, it is just
$\mu$ on each copy of $Y$.

Clearly, $(\widehat{Y},Y) $ is a ${\bf Z}/k$ manifold in the
sense of Sullivan (cf. [FM] and [Z1]).

Let $g^{\widehat{Y}}$ be a Riemannian metric on $T\widehat{Y}$
which is of product nature near $\partial \widehat{Y}$ and which
on $\partial \widehat{Y}$ is exactly $g^{TY}$ on each copy of the
boundary. Let $\nabla^{\widehat{Y}}$ be the associated
Levi-Civita connection. Similarly, let $g^{\widehat{\mu}}$ (resp.
$\nabla^{\widehat{\mu}}$) be a Hermitian metric on $\widehat{\mu}$
such that it is of product nature near $\partial \widehat{Y}$ and
that on $\partial \widehat{Y}$ it is exactly $g^\mu$ (resp.
$\nabla^{{\mu}}$)  on $\mu$ over each copy of $Y$.

Let $S^{2n,k}$ be the ${\bf Z}/k$ manifold obtained by removing
$k$ balls $D^{2n}$ from the $2n$-sphere. Then the boundary
$\partial S^{2n,k}$ consists of $k$ disjoint copies of $S^{2n-1}$.
Let $i:\widehat{Y}\hookrightarrow S^{2n,k}$ be a ${\bf Z}/k$
embedding (cf. [FM], [Z1]). The existence of such an embedding is
clear when $n$ is sufficiently large.

Let $g^{TS^{2n,k}}$ be a ${\bf Z}/k$ metric on $TS^{2n,k}$ which
is of product nature near $\partial S^{2n,k}$, such that
$g^{TS^{2n,k}}|_{T\widehat{Y}}=g^{T\widehat{Y}}$ and that the
isometric embedding
$i:(\widehat{Y},g^{T\widehat{Y}})\hookrightarrow(S^{2n,k},g^{TS^{2n,k}})$
is totally geodesic. Let $\nabla^{TS^{2n,k}}$ be the associated
Levi-Civita connection.

Let $\widehat{N}$ denote the normal bundle to $\widehat{Y}$ in
$S^{2n,k}$. Then $\widehat{N}$ carries an induced orientation and
spin structure, as well as a ${\bf Z}/k$ Euclidean metric (resp.
connection) $g^{\widehat{N}}$ (resp. $ \nabla^{\widehat{N}} $).

We can then apply the constructions in Sections 2a), b) to the
embedding
$i:(\widehat{Y},g^{T\widehat{Y}})\hookrightarrow(S^{2n,k},g^{TS^{2n,k}})$
in a ${\bf Z}/k$ manner (that is, preserving all the ${\bf Z}/k$
structures), such that all the metrics, connections and maps
involved are of product nature near boundary.

We denote the resulting ${\bf Z}/k$ geometric direct image of
$(\widehat{\mu},\nabla^{\widehat{\mu}})$ by
$(\widehat{\xi}=\widehat{\xi}_+\oplus\widehat{\xi}_-,\nabla^{\widehat{\xi}},\widehat{V})$.
When restrict to each copy of the boundary,  we denote  the
restricted geometric direct image of $(\mu,\nabla^\mu)$ by
$({\xi}={\xi}_+\oplus{\xi}_-,\nabla^{{\xi}},{V})$.

Let ${\rm ind}_k( D^{\widehat{\xi}_\pm})$, ${\rm ind}_k
(D^{\widehat{\mu}})$ be the mod $k$ indices defined by (cf. [FM,
(5.2)], [Z1]),
$${\rm ind}_k (D^{\widehat{\xi}_\pm}) \equiv\int_{S^{2n,k}}
\widehat{A}(TS^{2n,k},\nabla^{TS^{2n,k}}){\rm
ch}({\widehat{\xi}_\pm},\nabla^{\widehat{\xi}_\pm})-k\,\bar{\eta}(D^{{\xi}_\pm})
\ \ \ {\rm mod}\ \ k{\bf Z},\eqno(3.1)$$
$${\rm ind}_k (D^{\widehat{\mu}}) \equiv\int_{\widehat{Y}}
\widehat{A}(T\widehat{Y},\nabla^{\widehat{Y}}){\rm
ch}(\widehat{\mu},\nabla^{\widehat{\mu}})-k\,\bar{\eta}(D^{{\mu}})
\ \ \ {\rm mod}\ \ k{\bf Z}.\eqno(3.2)$$

By the mod $k$ index theorem of Freed and Melrose [FM, (5.5)], one
knows that
$${\rm ind}_k (D^{\widehat{\xi}_+})- {\rm ind}_k (D^{\widehat{\xi}_-})={\rm ind}_k
(D^{\widehat{\mu}})\eqno(3.3)$$ in ${\bf Z}/k{\bf Z}$.

Now if we denote $\gamma^{\widehat{\xi},\widehat{V}}$ the
Chern-Simons current constructed in Section 2b), then its
restriction to the boundary consists of $k$ copies of the
Chern-Simons current $\gamma^{\xi,V}$.

By applying the transgression formula (2.10) to
$\gamma^{\widehat{\xi},\widehat{V}}$ and integrate over
$S^{2n,k}$, one then gets
$$\int_{S^{2n,k}}
\widehat{A}(TS^{2n,k},\nabla^{TS^{2n,k}}){\rm
ch}({\widehat{\xi}_+},\nabla^{\widehat{\xi}_+})-\int_{S^{2n,k}}
\widehat{A}(TS^{2n,k},\nabla^{TS^{2n,k}}){\rm
ch}({\widehat{\xi}_-},\nabla^{\widehat{\xi}_-})$$
$$-\int_{\widehat{Y}}
\widehat{A}(T\widehat{Y},\nabla^{\widehat{Y}}){\rm
ch}(\widehat{\mu},\nabla^{\widehat{\mu}})
=k\int_{S^{2n-1}}\widehat{A}(TY,\nabla^{TY})\gamma^{\xi,V}.\eqno(3.4)$$

{}From (3.1)-(3.4), one deduces that
$$\bar{\eta}(D^{{\xi}_+})-\bar{\eta}(D^{{\xi}_-})\equiv\bar{\eta}
(D^{{\mu}})+\int_{S^{2n-1}}\widehat{A}(TY,\nabla^{TY})\gamma^{\xi,V}\
\ \ \ {\rm mod}\ \ {\bf Z} ,\eqno(3.5)$$ which is exactly the
Bismut-Zhang formula (2.11) in the case where $X=S^{2n-1}$.

$\ $

\noindent {\bf Remark 3.1} The relation between the Bismut-Zhang
formula (2.11) and the Freed-Melrose mod $k$ index theorem [FM]
was exploited in [Z1] where (2.11) is used to give an alternate
proof of a mod $k$ equality between the right hand sides of (3.1)
and (3.2). Now such an equality can be proved directly by
applying the Riemann-Roch property for Dirac operators on
manifolds with boundary proved by Dai and Zhang in [DZ], without
using the results in [BZ]. In fact, this can be done by first
applying [DZ, Theorem 1.2 and Lemma 4.6] to get (3.3). The mod
$k$ equality between the right hand sides of (3.1) and (3.2) is
then an easy consequence of the Atiyah-Patodi-Singer index
theorem [APS1] (This observation grew out of discussions with
Xianzhe Dai).

$\ $

\noindent {\bf Remark 3.2} One should also note that the proof of
(3.5) given in this section holds only for those embeddings which
can be obtained through  an embedding between ${\bf Z}/k$
manifolds, while (2.11) is much more general. Still we hope the
simplified proof in this section would be helpful for a good
feeling of (2.11). While on the other hand, for many
applications, the existence for such an embedding is suffice, as
our main result (2.19) holds for it.

$$\ $$

\noindent {\bf \S 4. Some applications}

$\ $

In this section, we discuss some immediate applications of
formula (2.19).

This section is organized as follows. In a), we discuss briefly
the relationship between the Chern-Simons current and the
Cheeger-Simons differential character [CS]. In b), we discuss the
Atiyah-Patodi-Singer ${\bf R}/{\bf Z}$ index theorem for unitary
flat vector bundles [APS3] from the point of view of (2.19). In
c), we show that by replacing the Dirac operator by the Signature
operator, one may refined the above ${\bf R}/{\bf Z}$ formula to
an ${\bf R }$ valued one.

We  make the same assumptions and use the same notation as in
Section 2.

$\ $

{\bf a). Chern-Simons current and the Cheeger-Simons differential
character }

$\ $

 As was indicated by
Bismut in [B2], the Chern-Simons currents constructed in Section
2 are closely related to the differential characters introduced
by Cheeger and Simons in [CS]. This becomes clearer if we compare
the transgression formula (2.18) with the one in [CS, (4.3)]. The
difference is that in [CS, (4.3)], the transgression formula
holds on different Stiefel manifolds, while our formula holds
universally on a single sphere. Moreover, the differential
characters for Chern character forms in [CS] were defined mod
${\bf Q}$, while our formula is clearly of an ${\bf R}/{\bf Z}$
nature (as the construction of the Chern-Simons current
$\gamma^{\widetilde{\xi},\widetilde{V}}$ in Section 2d) depends
on the choice of an automorphism $W$).

More precisely, if we denote by
$\widehat{\widehat{A}}(TY,\nabla^{TY})$ (resp. $\widehat{\rm
ch}(\mu,\nabla^\mu) $) the Cheeger-Simons differential character
associated to ${\widehat{A}}(TY,\nabla^{TY})$ (resp. ${\rm
ch}(\mu,\nabla^\mu) $) constructed in [CS], then by [CS, Theorem
9.1], one has, in using the product notation as in [CS],
$$\bar{\eta}(D^\mu)\equiv \left\langle
\widehat{\widehat{A}}(TY,\nabla^{TY}) *\widehat{\rm
ch}(\mu,\nabla^\mu) ,[Y]\right\rangle\ \ \ \ {\rm mod}\ \ {\bf
Q}.\eqno(4.1)$$

{}From (2.19) and (4.1), one gets
$$\int_{S^{2n-1}}\widehat{A}(TS^{2n-1},\nabla^{TS^{2n-1}})
\gamma^{\widetilde{\xi},\widetilde{V}}+ \left\langle
\widehat{\widehat{A}}(TY,\nabla^{TY}) *\widehat{\rm
ch}(\mu,\nabla^\mu) ,[Y]\right\rangle\equiv 0\ \ \ \ {\rm mod}\ \
{\bf Q}.\eqno(4.2)$$

It would be interesting to give a direct proof of (4.2).

$\ $

{\bf b). The Atiyah-Patodi-Singer index theorem for flat bundles
revisited}

$\ $

In this section, we replace the Hermitian vector bundle $\mu$ in
Section 2 by $\mu\otimes \rho$, where $\rho$ is a unitary flat
vector bundle with the flat connection denoted by $\nabla^\rho$.
We equip $\mu\otimes \rho$ with the induced tensor product
Hermitian metric as well as the tensor product Hermitian
connection $\nabla^{\mu\otimes \rho}=\nabla^{\mu}\otimes{\rm
Id}_\rho +{\rm Id}_\mu\otimes\nabla^\rho $.

By [APS3],
$$\bar{\eta}_{\mu,\rho}:=\bar{\eta}(D^{\mu\otimes\rho})-{\rm
rk}(\rho)\bar{\eta}(D^{\mu})\ \ \ \  {\rm mod} \ \ {\bf
Z}\eqno(4.3)$$ is a smooth invariant with respect to
$(g^{TM},g^\mu,\nabla^\mu)$. Moreover, [APS3, Theorem 5.3]
provides a topological interpretation for this invariant.

We now examine this invariant by using (2.19).

Let
$(\widetilde{\xi}_{\rho},\nabla^{\widetilde{\xi}_{\rho}},\widetilde{V}_{\rho})$
be a geometric direct image of $(\mu\otimes\rho,
\nabla^{\mu\otimes\rho})$ constructed similarly as that for
$(\mu,\nabla^\mu)$ in Section 2d). Let
$\gamma^{\widetilde{\xi}_{\rho},\widetilde{V}_{\rho}}$ be the
associated Chern-Simons current defined in (2.9). By (2.18),
$\gamma^{\widetilde{\xi}_{\rho},\widetilde{V}_{\rho}}$ verifies
the transgression formula
$$d\gamma^{\widetilde{\xi}_{\rho},\widetilde{V}_{\rho}}=-\widehat{A}^{-1}(N,\nabla^N)
{\rm ch}(\mu,\nabla^\mu){\rm rk}(\rho)\delta_Y,\eqno(4.4)$$ as
$\rho$ is a flat bundle.

By (2.19), one also has
$$\bar{\eta}(D^{\mu\otimes\rho})\equiv -\int_{S^{2n-1}}\widehat{A}(TS^{2n-1},\nabla^{TS^{2n-1}})
\gamma^{\widetilde{\xi}_{\rho},\widetilde{V}_{\rho}}\ \ \ \ {\rm
mod}\ \ {\bf Z}.\eqno(4.5)$$

{}From (2.19), (4.3) and (4.5), one gets
$$\bar{\eta}_{\mu,\rho}\equiv\int_{S^{2n-1}}\widehat{A}(TS^{2n-1},\nabla^{TS^{2n-1}})\left(
{\rm rk}(\rho)\gamma^{\widetilde{\xi},\widetilde{V}}-
\gamma^{\widetilde{\xi}_{\rho},\widetilde{V}_{\rho}}\right)\ \ \
\ {\rm mod}\ \ {\bf Z}.\eqno(4.6)$$

On the other hand, from (2.18) and (4.4), one finds
$$d\left(
{\rm rk}(\rho)\gamma^{\widetilde{\xi},\widetilde{V}}-
\gamma^{\widetilde{\xi}_{\rho},\widetilde{V}_{\rho}}\right)=0.\eqno(4.7)$$
Thus, ${\rm rk}(\rho)\gamma^{\widetilde{\xi},\widetilde{V}}-
\gamma^{\widetilde{\xi}_{\rho},\widetilde{V}_{\rho}}$ determines
a cohomology class in $H^{\rm odd}_{\rm dR}(S^{2n-1},{\bf R})$.
One verifies easily that this cohomology class does not depend on
the choice of $\nabla^{\mu}$. It only depends on the choices of
the automorphisms $W$ and $W_\rho$ appearing in the construction
of the geometric direct images of
$(\widetilde{\xi},\nabla^{\widetilde{\xi}},\widetilde{V})$ and
$(\widetilde{\xi}_{\rho},\nabla^{\widetilde{\xi}_{\rho}},\widetilde{V}_{\rho})$.
The different  choices of $W$ and $W_\rho$ cause a (possible)
integer jump in integration term in the right hand side of (4.6).

Thus, (4.6) may be thought of in some sense as an intrinsic
version of the Atiyah-Patodi-Singer index theorem for flat vector
bundles stated in [APS3, Theorem 5.3]. Its conceptual novelty is
that one need not divide the topological index into two parts
(that is, a ${\bf Q}/{\bf Z}$ part plus an ${\bf R}$ part.)

$\ $

\noindent {\bf Remark 4.1} It is pointed out in [APS3] that the
index theorem for flat vector bundles stated in [APS3, Theorem
5.3] indeed holds also for non-unitary flat vector bundles by
suitable extension of the construction of the $\eta$ invariants.
In the next section, we will extend our intrinsic formula (4.6)
to the case of non-unitary flat vector bundles along similar
lines.

$\ $

{\bf c). Signature operator and an ${\bf R}$ valued index theorem
for flat vector bundles }

$\ $

Now we set $\mu=S(TY)$ in the above subsection. In this case,
$D^\mu$ is the Signature operator associated to $(TY,g^{TY})$,
denoted by $D_{\rm Sign}$, while $D^{\mu\otimes\rho}$ is now
denoted by $D^\rho_{\rm Sign}$.

Set
$$\bar{\eta}_{{\rm Sign},\rho}=\bar{\eta}(D^\rho_{\rm Sign})-{\rm rk}( \rho)
\bar{\eta}(D_{\rm Sign}).\eqno(4.8)$$ Then $\bar{\eta}_{{\rm
Sign},\rho}$ is a smooth invariant equivalent to what defined in
[APS2, Theorem 2.4].

On the other hand, as indicated in Remark 2.4, one can choose the
automorphisms $W$ and $W_\rho$ in the construction of the
Chern-Simons current so that
$$\bar{\eta}(D^\rho_{\rm Sign})=\int_{S^{2n-1}}\widehat{A}(TS^{2n-1},\nabla^{TS^{2n-1}})
\gamma^{\widetilde{\xi}_{\rho},\widetilde{V}_{\rho}},\eqno(4.9)$$
$$\bar{\eta}(D_{\rm
Sign})=\int_{S^{2n-1}}\widehat{A}(TS^{2n-1},\nabla^{TS^{2n-1}})
\gamma^{\widetilde{\xi},\widetilde{V}}.\eqno(4.10)$$ It is clear
that (4.9), (4.10) are actually equalities not depending on the
choice of $g^{TM}$. Moreover, since $K^1(S^{2n-1})={\bf Z}$, one
sees that the choice of $W$ and $W_\rho$ is canonical (up to
stable homotopy).

By (4.7), one also sees that
$\gamma^{\widetilde{\xi}_{\rho},\widetilde{V}_{\rho}}-{\rm rk
}(\rho)\gamma^{\widetilde{\xi},\widetilde{V}}$ determines
canonically an element in $H^{\rm odd}_{\rm dR}(Y,{\bf R})$.

We can state our ${\bf R}$ valued refinement of (4.6) 
 as follows.

$\ $

\noindent {\bf Theorem 4.2} {\it The following $K$-theoretic
formula for $\bar{\eta}_{{\rm Sign},\rho}$ holds,
$$\bar{\eta}_{{\rm Sign},\rho}=\int_{S^{2n-1}}\widehat{A}(TS^{2n-1},\nabla^{TS^{2n-1}})\left(
{\rm rk}(\rho)\gamma^{\widetilde{\xi},\widetilde{V}}-
\gamma^{\widetilde{\xi}_{\rho},\widetilde{V}_{\rho}}\right).\eqno(4.11)$$}

\noindent {\bf Remark 4.3} Theorem 4.2 was proved for unitary flat
vector bundles over spin manifolds. It is nature to ask whether
there is  still such a kind of formulas without the spin
condition.

$\ $

\noindent {\bf Remark 4.4} If one could find a purely topological
way to identify the automorphisms $W$ and $W_\rho$ (or  the
difference element of ${\rm rk}(\rho)W$ and $W_\rho$ in
$K^1(S^{2n-1})$), then one would provide a positive answer to a
question of Atiyah-Patodi-Singer stated implicitly in [APS2, page
406]. On the other hand, for any choice of $W$ and $W_\rho$, the
right hand side of (4.11) provides a smooth invariant of $Y$. So
in some sense, $\bar{\eta}_{{\rm Sign},\rho}$ becomes one example
of a series of smooth invariants associated to the unitary flat
vector bundle $\rho$ over $Y$.

$$\ $$

\noindent {\bf\S 5. Generalizations to Hermitian vector
 bundles with non-unitary connection}

$\ $

In this section, we generalize the results in the previous
sections to the case of non-unitary (i.e., non-Hermitian)
connections on the twisted Hermitian vector bundle $\mu$.

This section is organized as follows. In a), we propose a
definition of the $\eta$ invariants of Dirac operators coupled
with Hermitian vector bundles with non-unitary connection. In b),
we show that the generalized $\eta$ invaraints defined in a) can
also be expressed as a Chern-Simons current on a higher
dimensional sphere, extending Theorem 2.3. In c), we apply the
formula obtained in b) to the case of general flat vector bundles.

$\ $

{\bf a). $\eta$ invariants for Hermitian vector
 bundles with non-unitary connection}

 $\ $

 We make the same assumptions and use the same notation as in
 Sections 1 and 2, except that we no longer assume that the connection
 $\nabla^\mu$ is a Hermitian connection with respect to the
 Hermitian metric $g^\mu$. Then the Dirac operator $D^\mu$ need
 not be self-adjoint. To make clearer the dependence of the Dirac
 operator with respect to the connection, we will denote the Dirac
 operator coupled with connection $\nabla^\mu$ on $\mu$ by
 $D^{\nabla^\mu}$.

 Let $\nabla^{\mu,*}$ be the adjoint connection of $\nabla^\mu$
 with respect to $g^\mu$. Then
 $$\widetilde{\nabla}^\mu={1\over
 2}\left(\nabla^\mu+\nabla^{\mu,*}\right)\eqno(5.1)$$
 is a Hermitian connection on $\mu$, and
 $D^{\widetilde{\nabla}^\mu}$ is self-adjoint.

 Let $\psi :{\bf R}\rightarrow {\bf R}$ be a smooth function such
 that $\psi(x)=1$ if $x>10$, while $\psi(x)=0$ if $x<2$.

 For any $t\geq 0$, set
 $$D(t)=\sqrt{t}\left(\psi(t) D^{\widetilde{\nabla}^\mu}+(1-\psi(t)
 )D^{\nabla^\mu}\right).\eqno(5.2)$$

 First of all, by using Duhamel principle, one knows easily that
 for any $t>0$, the heat operator $\exp(-D(t)^2)$ is well-defined.

 By proceeding similarly as in [BF, Theorem 2.4], one knows that as
 $t\rightarrow 0^+$,
 $${\rm Tr}\left[{dD(t)\over
 dt}\exp\left(-D(t)^2\right)\right]=O(1).\eqno(5.3)$$

Thus, one is able to define, by analogous to [L, (59)],
$${\eta}(D^{\nabla^\mu},\psi)={2\over\sqrt{\pi}}\int_{0}^{+\infty}
{\rm Tr}\left[{dD(t)\over
 dt}\exp\left(-D(t)^2\right)\right].\eqno(5.4)$$

 By proceeding as in [L, Lemma 3], one knows that
 ${\eta}(D^{\nabla^\mu},\psi)$ does not depend on $\psi$.

 $\ $

 \noindent {\bf Definition 5.1} {Let ${\eta}(D^{\nabla^\mu})$
 denote ${\eta}(D^{\nabla^\mu},\psi)$ and call it the $\eta$
 invariant of $D^{\nabla^\mu}$. Let $\bar{\eta}(D^{\nabla^\mu})$
 be the reduced $\eta$ invariant defined by
 $$\bar{\eta}(D^{\nabla^\mu})={\dim(\ker D^{\widetilde{\nabla}^\mu})+{\eta}(D^{\nabla^\mu})\over
 2}.\eqno(5.5)$$}

 Clearly, if $\nabla^\mu$ preserves $g^\mu$, then
 $\bar{\eta}(D^{\nabla^\mu})$ is exactly the usual reduced $\eta$
 invariant discussed in Sections 1 and 2.

 Let $CS(\nabla^\mu,\widetilde{\nabla}^\mu)$ be the
 Chern-Simons form verifying (cf. [Z2, Chapter 1])
 $$dCS(\nabla^\mu,\widetilde{\nabla}^\mu)={\rm
 ch}(\mu,\nabla^\mu)-{\rm
 ch}(\mu,\widetilde{\nabla}^\mu).\eqno(5.6)$$

 $\ $

 \noindent {\bf Proposition 5.2} {\it The following identity
 holds,}
 $$\bar{\eta}(D^{\nabla^\mu})=\bar{\eta}(D^{\widetilde{\nabla}^\mu})+\int_Y
 \widehat{A}(TM,\nabla^{TM})CS(\nabla^\mu,\widetilde{\nabla}^\mu).\eqno(5.7)$$

 {\it Proof}. For any $u\in [0,1]$, set $D^{\nabla^\mu}(u)=(1-u)D^{\nabla^\mu}+uD^{\widetilde{\nabla}^\mu}$.
 By proceeding as in [BF, Theorem 2.10], one knows that the ${\bf
 R}/{\bf Z}$ part of
 $\bar{\eta}(D^{\nabla^\mu}(u))$ verifies the variation formula as usual.
 This gives the local term in the right hand side of (5.7).
 On the other hand, it is clear that there is no integer jump in
 this deformation (Compare with [L, (65)]). \ \ Q.E.D.

$\ $

{\bf b). A formula through Chern-Simons current}

$\ $

In this subsection, we show that Theorem 2.3 still holds in the
current generalized situation. We state the result as follows.

$\ $

\noindent {\bf Theorem 5.3} {\it The following identity holds,}
$$\bar{\eta}(D^\mu)\equiv -\int_{S^{2n-1}}\widehat{A}(TS^{2n-1},\nabla^{TS^{2n-1}})
\gamma^{\widetilde{\xi},\widetilde{V}}\ \ \ {\rm mod}\ {\bf
Z}.\eqno(5.8)$$

{\it Proof}. Let
$(\widetilde{\xi},\widehat{\nabla}^{\widetilde{\xi}},\widetilde{V}
)$ be the geometric direct image constructed in Section 2d) with
respect the connection $\widetilde{\nabla}^\mu$ on $\mu$ instead
of $\nabla^\mu$. Let
$\widehat{\gamma}^{\widetilde{\xi},\widetilde{V}}$ be the
associated Chern-Simons current. Then by Theorem 2.3, we know that
$$\bar{\eta}(D^{\widetilde{\nabla}^\mu})\equiv -\int_{S^{2n-1}}\widehat{A}(TS^{2n-1},\nabla^{TS^{2n-1}})
\widehat{\gamma}^{\widetilde{\xi},\widetilde{V}}\ \ \ {\rm mod}\
{\bf Z}.\eqno(5.9)$$

On the other hand, by the transgression formula (2.18), one
deduces that
$$\int_{S^{2n-1}}\widehat{A}(TS^{2n-1},\nabla^{TS^{2n-1}})
\widehat{\gamma}^{\widetilde{\xi},\widetilde{V}}-\int_{S^{2n-1}}\widehat{A}(TS^{2n-1},\nabla^{TS^{2n-1}})
{\gamma}^{\widetilde{\xi},\widetilde{V}}$$ $$=\int_Y
 \widehat{A}(TM,\nabla^{TM})CS(\nabla^\mu,\widetilde{\nabla}^\mu).\eqno(5.10)$$

 {}From (5.7), (5.9) and (5.10), one gets (5.8).\ \ Q.E.D.

 $\ $

 {\bf c). Applications to flat vector bundles}

 $\ $

 We will only be brief in this subsection as the arguments run
 parallel to those in Section 4.

 Now we replace $\mu$  by $\mu\otimes\rho$ such that $\mu$ is a
 Hermitian vector bundle carrying with a Hermitian connection,
 while $\rho$ is a flat vector bundle carrying with a Hermitian
 metric. But we no longer assume that the flat connection
 $\nabla^\rho$ is unitary.

 With the definition given in Section 5a), it is clear that the
 invariant $\bar{\eta}_{\mu,\rho}$ is still well-defined.
 Moreover, by Theorem 5.3,
 one still has a formula of the same form as (4.6), giving a ${\bf
 C}/{\bf Z}$ valued
 $K$-theoretic formula for this analytically defined invariant. It is
 easy to check that the right hand side of such a formula is equal to
  the topological index given  in [APS3, Theorem
 5.3]. This means that, at least on the mod ${\bf Z}$ level, our definition of the
 reduced $\eta$ invariant agrees with the one indicated in [APS3, Sections
 5, 6]. This would be sufficient for  many other applications as well.

 Certainly one can also extend Theorem 4.2 to the case of
 non-unitary flat vector bundles.

 We leave the details  to the
 interested reader.

 $\ $

 \noindent {\bf Remark 5.4} Clearly, all the results of this
 paper can be extended to the case of spin$^c$ manifolds. We also leave
 this to the interested reader.

$$\ $$

\noindent {\bf Acknowledgement} The author is indebted to
Professor Jean-Michel Bismut in an obvious way. He would also like
to thank Fuquan Fang, Huitao Feng and Xiaonan Ma for helpful
discussions.

$$\ $$

\noindent {\bf References}

$\ $

\noindent [A] M. F. Atiyah, {\it K-theory}. Benjamin, New York,
1967.

$\ $

\noindent [AH] M. F. Atiyah and F. Hirzebruch, Riemann-Roch
theorems for differentiable manifolds. {\it Bull. Amer. Math.
Soc.} 65 (1959), 276-281.

$\ $

\noindent [APS1] M. F. Atiyah, V. K. Patodi and I. M. Singer,
Spectral asymmetry and Riemannian geometry I. {\it Proc.
Cambridge Philos. Soc.} 77 (1975), 43-69.

$\ $

\noindent [APS1] M. F. Atiyah, V. K. Patodi and I. M. Singer,
Spectral asymmetry and Riemannian geometry II. {\it Proc.
Cambridge Philos. Soc.} 78 (1975), 405-432.

 $\ $

\noindent [APS3] M. F. Atiyah, V. K. Patodi and I. M. Singer,
Spectral asymmetry and Riemannian geometry III. {\it Proc.
Cambridge Philos. Soc.} 79 (1976), 71-99.

$\ $

\noindent [AS] M. F. Atiyah and I. M. Singer, The index of
elliptic operators I. {\it Ann. of Math.} 87 (1968), 484-530.

$\ $

\noindent [B1] J.-M. Bismut, Superconnection currents and complex
immersions. {\it Invent. Math.} 99 (1990), 59-113.

$\ $

\noindent [B2] J.-M. Bismut, Eta invariants and complex
immersions. Bull. Soc. Math. France 118 (1990), 211-227.

$\ $

\noindent [BC] J.-M. Bismut and J. Cheeger, $\eta$-invariants and
their adiabatic limits. {\it J. Amer. Math. Soc.} 2 (1989), 33-70.

$\ $

\noindent [BF] J.-M. Bismut and D. S. Freed, The analysis of
elliptic families, II. {\it Commun. Math. Phys.} 107 (1986),
103-163.

$\ $

 \noindent [BGS] J.-M. Bismut, H. Gillet and C. Soul\'e,
Bott-Chern currents and complex geometry. {\it Duke Math. J.} 60
(1990), 255-284.

$\ $

\noindent [BL] J.-M. Bismut  and G. Lebeau, Complex immersions
and Quillen metrics, {\em Publ. Math. IHES.} 74 (1991), 1-297.

$\ $

\noindent [BZ] J.-M. Bismut and W. Zhang, Real embeddings and eta
invariant. {\it Math. Ann.} 295 (1993), 661-684.

$\ $

\noindent [CS] J. Cheeger and J. Simons, Differential characters
and geometric invariants. in {\it Lecture Notes in Math.} Vol.
1167, page 50-80. Springer-Verlag, 1985.

$\ $

\noindent [DZ] X. Dai and W. Zhang,  Real embeddings and the
Atiyah-Patodi-Singer index theorem for Dirac operators. {\it Asian
J. Math.} 4 (2000), 775-794.

$\ $

\noindent [FM] D. S. Freed and R. B. Melrose, A mod $k$ index
theorem. {\it Invent. Math.} 107 (1992), 283-299.

$\ $

\noindent [L] J. Lott, ${\bf R}/{\bf Z}$-index theory. {\it
Commun. Anal. Geom.} 2 (1994), 279-311.

$\ $

\noindent [Q] D. Quillen, Superconnections and the Chern
character. {\it Topology} 24 (1985), 89-95.

$\ $

\noindent [Z1] W. Zhang, On the mod $k$ index theorem of Freed and
Melrose. {\it J. Diff. Geom.} 43 (1996), 198-206.

$\ $

\noindent [Z2] W. Zhang, {\it Lectures on Chern-Weil Theory and
Witten Deformations.} Nankai Tracks in Mathematics Vol. 4, World
Scientific, Singapore, 2001.

$$\ $$

Nankai Institute of Mathematics,  Nankai University, Tianjin
300071,  P. R. China

 {\it E-mail: weiping@nankai.edu.cn}
\end{document}